\documentclass{article}
\usepackage[utf8]{inputenc}
\usepackage[english]{babel}
\usepackage{amsmath}
\usepackage{amsthm}
\usepackage{amsfonts}
\usepackage{amssymb}
\usepackage{enumitem}

\begin{document}

\author{Walter Wyss}
\title{Heron Angles, Heron Triangles, and Heron Parallelograms} 
\date{}
\maketitle

\begin{abstract}
\noindent Heron angle: both its sine and cosine are rational \\
Heron triangle: all its sides and area are rational \\
Heron Parallelogram: all its sides, diagonals and area are rational \\
We give one-to-one (bijective) parametrizations for all three concepts.
\end{abstract}

\section{Introduction}
Special angles, triangles, rectangles, and parallelograms have been studied for ages by Mathematician and Scientists (Pythagoras, Heron, Diophantus, Brahmagupta to mention a few). Parameter representations, some in integers, have been found in all cases \cite{1}. However, there was a lack in retrieving the parameters. Here we give parameter representations that are bijective (one-to-one). 

\section{Bijective parameter representation of a relation}

\subsection*{\underline{Theorem 1}} let $a,x,y,z$ be real numbers. Given $a$, the relation

\begin{equation}
x^2 + 2axy + y^2 = z^2
\end{equation}
 has the following bijective parameter representation, with parameters $ \lambda, \sigma$
 
\begin{align}
x &= \frac{\lambda}{1+a\sigma}(1-\sigma^2) \\
y &= 2\lambda\sigma \\
z &= \frac{\lambda}{1+a\sigma}(1+2a\sigma+\sigma^2) \\
\text	{conversely}  \hspace{0.5cm} 2\lambda &= x+z, \hspace{0.5cm} \sigma = \frac{y}{x+z}
\end{align}

\begin{proof}
The relation (1) also reads
\begin{equation}
x^2 + (y+ax)^2 = z^2 + (ax)^2
\end{equation}
\end{proof}

According to \cite{2} this relation has the following bijective parameter representation with parameters $s_1, s_2, \lambda_2$

\begin{align}
y+ax = s_1+\lambda_2s_2 \hspace{1.0cm} ax = s_1-\lambda_2s_2 \\
x= s_2-\lambda_2s_1 \hspace{1.0cm} z=s_2+\lambda_2s_1
\end{align}

Conversely
\begin{align}
2s_1 &= y+2ax \\
2s_2 &= x+z \\
2\lambda_2s_2 &= y
\end{align}

Now
\begin{align}
ax = a_1(s_2-\lambda_2s_1) = s_1-\lambda_2s_2
\end{align}

resulting in
\begin{align}
s_1 = \frac{a+\lambda_2}{1+a\lambda_2}s_2
\end{align}

Then
\begin{align}
\begin{split}
x &= s_2 - \lambda_2 \frac{a+\lambda_2}{1+a\lambda_2}s_2 \\
x &= \frac{s_2}{1+a\lambda_2} (1-\lambda_2^2) 
\end{split}
\end{align}
\begin{align}
\begin{split}
z &= s_2 + \lambda_2 \frac{a+\lambda_2}{1+a\lambda_2} s_2 \\
z &= \frac{s_2}{1+a\lambda_2} (1+2a\lambda_2 + \lambda_2^2)
\end{split}
\end{align}

and
\begin{align}
y &= 2\lambda_2s_2 
\end{align}

Finally let  $s_2 = \lambda, \hspace{0.15cm} \lambda_2 = \sigma$

\subsection*{\underline{Corollary 1}} The relation, given $\varepsilon$, 
\begin{equation}
m^2 + \varepsilon mn + n^2 = 1
\end{equation}

has the following bijective parameter representation, with parameter $\sigma$

\begin{align}
m &= \frac{1-\sigma^2}{1+\varepsilon\sigma + \sigma^2} \\
n &= \frac{\sigma (2+\varepsilon \sigma)}{1 + \varepsilon\sigma+\sigma^2}
\end{align}

Conversely
\begin{align}
\sigma &= \frac{n}{1+m}
\end{align}

\begin{proof}
In (2, 3, 4, 5) let $x =m, \hspace{0.15cm} y =n, \hspace{0.15cm} z=1, $ and $2a = \varepsilon$
\end{proof}

\section{The relation}
\begin{equation}
\lambda^2 = \frac{m}{n} \frac{1-n^2}{1+m^2}
\end{equation}

We have two one-parameter families of fundamental solutions given by

\begin{enumerate}[label=(\alph*)]
\item Type I: \begin{equation}
\lambda = \frac{m}{n}
\end{equation}
\newline
Then
\begin{align*}
m(1-m^2) &= n(1-n^2) \\
m-n &= m^3 - n^3
\end{align*}

or
\begin{equation}
m^2 +mn +n^2 = 1
\end{equation}

This relation has already been considered by Diophantus of Alexandria \cite{3} \\
According to (17), the relation (23) has the bijective parameter representation with parameter $\sigma$, as 

\begin{equation}
m=\frac{1-\sigma^2}{1+\sigma+\sigma^2}, \hspace{0.5cm} n=\frac{\sigma(2+\sigma)}{1+\sigma+\sigma^2}, \hspace{0.5cm} \lambda=\frac{1-\sigma^2}{\sigma(2+\sigma)}
\end{equation}

conversely
\begin{equation}
\sigma=\frac{n}{1+m}
\end{equation}

\item Type II: \begin{equation}
\lambda = mn
\end{equation}
\newline
Then
\begin{align*}
mn^3 (1-m^2) &= 1-n^2 \\
1+m^3n^3 &= n^2(1+mn)
\end{align*}
or
\begin{equation*}
1-mn+m^2n^2=n^2
\end{equation*}
and then 
\begin{equation}
\bigg(\frac{1}{n}\bigg)^2 - m\bigg(\frac{1}{n}\bigg) +m^2 =1
\end{equation}
\end{enumerate}

According to (17), the relation (27) has the bijective parameter representation with parameter $\sigma$, as

 \begin{equation}
m=\frac{1-\sigma^2}{1-\sigma+\sigma^2}, \hspace{0.5cm} n=\frac{1-\sigma+\sigma^2}{\sigma(2-\sigma)}, \hspace{0.5cm} \lambda=\frac{1-\sigma^2}{\sigma(2-\sigma)}
\end{equation}

conversely
\begin{equation}
\sigma=\frac{1}{n(1+m)}
\end{equation}

\subsection*{\underline{Example}}

\begin{enumerate}[label=(\alph*)]
\item
\begin{equation}
\sigma=\frac{1}{2}, \hspace{0.5cm} m=\frac{3}{7}, \hspace{0.5cm} n=\frac{5}{7}, \hspace{0.5cm} \lambda=\frac{3}{5}
\end{equation}
\item
\begin{equation}
\sigma=\frac{2}{3}, \hspace{0.5cm} m=\frac{5}{7}, \hspace{0.5cm} n=\frac{7}{8}, \hspace{0.5cm} \lambda=\frac{5}{8}
\end{equation}
\end{enumerate}

\section{Heron angles}

\subsection*{\underline{Definition 1}}

\begin{enumerate}
\item An angle $\alpha$ is called a Heron angle if both $\sin\alpha$ and $\cos\alpha$ are rational
\item The generator of an angle $\alpha$ is defined by
\begin{equation}
m = m(\alpha) = \frac{\sin \alpha}{1+\cos \alpha} = \tan\bigg(\frac{\alpha}{2}\bigg)
\end{equation}
\end{enumerate}

\subsection*{\underline{Lemma 1}}

For $0 < \alpha < \pi, \hspace{0.20cm} m(\alpha) > 0$ and is an increasing function of $\alpha$. 

\begin{proof}
Since $\sin \alpha > 0$ and $-1 < \cos\alpha < 1$ we see from (32) that $m(\alpha) > 0 $. \\
\newline
From
\begin{equation*}
\frac{dm(\alpha)}{d\alpha} = \frac{1}{2} \frac{1}{\cos^2(\frac{\alpha}{2})} > 0
\end{equation*}
we see that $m(\alpha)$ is increasing.
\end{proof}

\subsection*{\underline{Lemma 2}}

Given the generator $m=m(\alpha), \hspace{0.2cm} 0 < \alpha < \pi$,
we find

\begin{equation}
\cos \alpha = \frac{1-m^2}{1+m^2}, \hspace{0.5cm} \sin \alpha = \frac{2m}{1+m^2}
\end{equation}

\begin{proof}
From 
\begin{equation*}
m^2 = \frac{1-\cos \alpha}{1+\cos \alpha}
\end{equation*}
we find
\begin{equation*}
\cos \alpha = \frac{1-m^2}{1+m^2} \hspace{0.5cm} \text{and then} \hspace{0.5cm} \sin \alpha = \frac{2m}{1+m^2}
\end{equation*}

Observe that the generator of a Heron angle is rational and vice versa. Thus there is a one-to-one relationship between rational numbers and Heron angles.

\end{proof}

\section{Heron triangles}

For a triangle with sides $u_1, u_2, u_3$ and interior angles $\phi_1, \phi_2, \phi_3$, where $\phi_k$ is the angle opposite $u_k$, $k=1,2,3$ \\
\newline
We have 
\begin{enumerate}
\item The law of sine
\begin{equation}
\frac{u_1}{\sin \phi_1} = \frac{u_2}{\sin \phi_2} + \frac{u_3}{\sin \phi_3}
\end{equation}
\item The law of cosine
\begin{equation}
u_3^2 = u_1^2 + u_2^2 - 2u_1u_2 \cos \phi_3
\end{equation}
with similar relations involving the other sides and angles
\item The area $A$ of the triangle is given by
\begin{equation}
A = \frac{1}{2}u_1u_2\sin \phi_3
\end{equation}
and similarly for the other sides and angles.
\end{enumerate}

\subsection*{\underline{Definition 2}}

A triangle is called a Heron triangle if all its sides and area are rational. \\
Observe that from (35, 36) all the interior angles are Heron angles. \\
For a Heron triangle we find from (34) and $\phi_1 = \pi - (\phi_2 +\phi_3)$ the representation

\begin{align}
u_1 &= \sin (\phi_2 + \phi_3) w \\
u_2 &= \sin \phi_2 \cdot w \\
u_3 &= \sin  \phi_3 \cdot w
\end{align}

where $w$ is a rational scaling parameter. \\
Let now $q$ be the generator of $\phi_2$, $p$ the generator of $\phi_3$, i.e.

\begin{align}
\sin \phi_2 = \frac{2q}{1+q^2}, \hspace{0.5cm} \cos \phi_2 = \frac{1-q^2}{1+q^2} \\
\sin \phi_3 = \frac{2p}{1+p^2}, \hspace{0.5cm} \cos \phi_3 = \frac{1-p^2}{1+p^2} 
\end{align}

Then
\begin{align*}
\sin \phi_1 &= \sin (\phi_2 + \phi_3) = \sin \phi_2 \cos \phi_3 + \cos \phi_2 \sin \phi_3 \\
\sin \phi_1 &= \frac{2q(1-p^2)}{(1+q^2)(1+p^2)} + \frac{2p(1-q^2)}{(1+q^2)(1+p^2)}
\end{align*}

\begin{equation}
\sin \phi_1 = \frac{2(p+q)(1-pq)}{(1+q^2)(1+p^2)}
\end{equation}

Introduce the new rational scaling parameter $v$ by

\begin{equation}
w = \frac{1}{2} (1+p^2)(1+q^2)v
\end{equation}

Then we have the representation, with $p>0, q>0, v>0 $ and $ pq<1$

\begin{align}
u_1 &= (p+q)(1-pq)v \\
u_2 &= q (1+p^2) v \\
u_3 &= p(1+q^2) v \\
A &= pq(p+q)(1-pq)v^2
\end{align}

This representation has the three parameters $p, q, v$. We now can retrieve these parameters from $u_1, u_2, u_3, A$ as follows:

From
\begin{align}
\sin \phi_3 &= \frac{4A}{2u_1u_2} \\
\cos \phi_3 &= \frac{u_1^2+u_2^2-u_3^2}{2u_1u_2} \\
1 + \cos \phi_3 &= \frac{(u_1+u_2)^2 - u_3^2}{2u_1u_2}
\end{align}

we find the generator
\begin{equation}
p = \frac{\sin \phi_3}{1+\cos \phi_3} = \frac{4A}{(u_1+u_2)^2-u_3^2}
\end{equation}

and similarly
\begin{equation}
q = \frac{\sin \phi_2}{1+\cos \phi_2} = \frac{4A}{(u_1+u_3)^2-u_2^2}
\end{equation}

The scaling factor $v$ is then given by
\begin{equation}
v=\frac{u_3}{p(1+q^2}
\end{equation}

This is our one-to-one relationship.

\subsection*{\underline{Example}}

\begin{align*}
u_1 = 9, \hspace{0.25cm} u_2 = 10, \hspace{0.25cm} u_3 = 17, \hspace{0.25cm} A=36 \\
p =2, \hspace{0.25cm} q =\frac{1}{4}, \hspace{0.25cm} v=8
\end{align*}

\section{Heron parallelograms}

A parallelogram has sides $u_1, u_2$ and diagonals  $u_3, u_4$.

\subsection*{\underline{Definition 3}}

\begin{enumerate}[label=(\alph*)]
\item A parallelogram with its sides and diagonals being rational is called a rational parallelogram.
\item A rational parallelogram with rational area is called a Heron parallelogram.
\end{enumerate}

In \cite{4} we found a bijective parameter representation for rational parallelograms. 

Parameters $u > 0, \hspace{0.25cm} 0 < m< 1,  \hspace{0.25cm} 0 < n <1$ and all rational

\begin{align}
u_1 &= (1-mn)u \\
u_2 &= (m+n)u \\
u_3 &= [ 1+mn - (n-m) ]u \\
u_4 &= [ 1+mn + (n-m) ] u
\end{align}

conversly
\begin{align}
4u &= u_4 +u_3 +2u_1 \\
m &= \frac{u_4 +u_3 -2u_1}{u_4 -u_3 +2u_2} \\
n &= \frac{u_4 -u_3 +2u_2}{u_4 +u_3 +2u_1}
\end{align}

Not let $\phi$ be then angle between the sides $u_1, u_2$.
According to the law of cosine we find

\begin{equation}
\cos \phi = \frac{u_1^2+u_2^2-u_3^2}{2u_1u_2} = \frac{2(u_1^2+u_2^2)-2u_3^2}{4u_1u_2}
\end{equation}

From the parallelogram equation
\begin{equation}
2(u_1^2+u_2^2) = u_3^2 +u_4^2
\end{equation}

we get
\begin{equation*}
\cos \phi = \frac{u_4^2-u_3^2}{4u_1u_2} = \frac{(u_4-u_3)(u_4+u_3)}{4u_1u_2}
\end{equation*}

\begin{equation}
\cos \phi = \frac{(n-m)(1+mn)}{(n+m)(1-mn)}
\end{equation}

Now the area $A$ of a parallelogram is given by
\begin{equation}
A=u_1u_2\sin \phi
\end{equation}

where 
\begin{align*}
\sin^2 \phi &=  1 - \cos^2 \phi \\
\sin^2 \phi &=  \frac{4mn(1-m^2)(1-n^2)}{(n+m)^2(1-mn)^2}
\end{align*}

\begin{align}
\sin^2 \phi &=  \frac{4n^2(1-m^2)^2}{(n+m)^2(1-mn)^2} \cdot \frac{m(1-n^2)}{n(1-m^2)}
\end{align}

or according to (21) with 
\begin{equation*}
\lambda^2 = \frac{m(1-n^2)}{n(1-m^2)}
\end{equation*}

we have 
\begin{equation}
\sin \phi = \frac{2n(1-m^2)}{(n+m)(1-mn)} \lambda
\end{equation}

The generator $p$ of the angle $\phi$ is given by
\begin{equation*}
p=\frac{\sin \phi}{1+\cos \phi}
\end{equation*}

where
\begin{equation*}
1 + \cos \phi = \frac{(n+m)(1-mn)+(n-m)(1+mn)}{(n+m)(1-mn)}
\end{equation*}
\begin{equation}
1 + \cos \phi = \frac{2n(1-m^2)}{(n+m)(1-mn)}
\end{equation}

Therefore 
\begin{equation}
p=\lambda
\end{equation}

The area is now given by
\begin{equation}
A=2\lambda n(1-m^2)u^2
\end{equation}

\noindent Now, a Heron parallelogram with sides $u_1, u_2$, diagonals $u_3, u_4$ and area $A$ is parameterized by
\begin{equation}
0<m<1, \hspace{0.25cm} 0<n<1, \hspace{0.25cm} u>0, \hspace{0.25cm} \lambda>0 \text{ with } \lambda^2=\frac{m(1-n^2)}{n(1-m^2)}
\end{equation}

\noindent For a type I Heron parallelogram, we have the relation (22) and thus the bijective parameter $\sigma$. From (24, 25, 70) we find

\begin{equation}
0<\sigma<1
\end{equation}

\noindent For a type II Heron parallelogram, we have the relation (26) and thus the bijective parameter $\sigma$. From (28, 29, 70) we find 

\begin{equation}
\frac{1}{2}<\sigma<1
\end{equation}

\noindent Observe that Heron parallelograms cover the case of Heron triangles with a rational median. \\
\newline
\noindent Finally, we have the relations 
\newline
\begin{flalign}
\begin{split}
\text{Type I :}  & \\ 
&\lambda = \frac{m}{n} = \frac{2u_2+u_3-u_4}{2u_2-u_3+u_4} \\
\text{Type II :} & \\ 
&\lambda = mn = \frac{u_3+u_4-2u_1}{u_3+u_4+2u_1}
\end{split}
\end{flalign}

\subsection*{\underline{Example}}

Type I :
From (30)
\begin{align}
\begin{split}
m &=\frac{3}{7}, \hspace{0.25cm} n=\frac{5}{7}, \hspace{0.25cm} \lambda=\frac{3}{5}, \hspace{0.25cm} u=\frac{49}{2} \\
u_1 &=17, \hspace{0.25cm} u_2=28, \hspace{0.25cm} u_3=25, \hspace{0.25cm} u_4=39, \hspace{0.25cm} A=420
\end{split}
\end{align}
\newline
Type II :
From (31)
\begin{align}
\begin{split}
m &=\frac{5}{7}, \hspace{0.25cm} n=\frac{7}{8}, \hspace{0.25cm} \lambda=\frac{5}{8}, \hspace{0.25cm} u=56 \\
u_1 &=21, \hspace{0.25cm} u_2=89, \hspace{0.25cm} u_3=82, \hspace{0.25cm} u_4=100, \hspace{0.25cm} A=1680
\end{split}
\end{align}

\bigskip

\noindent\textit{Department of Physics, University of Colorado Boulder, Boulder, CO 80309\\
Walter.Wyss@Colorado.EDU}


\bigskip

\begin{thebibliography}{100}
\bibitem{1} Hermann C.H. Schubert, {\em Die Ganzzahligkeit in der Algebrauschen Geometrie (1905)},
  Translated by Ralph H. Buchholz, Integrability in algebraic geometry (2005)
\bibitem{2} Walter Wyss, {\em Sum of Squares, Bijective Parameter Representation},\\\texttt{https://arxiv.org/abs/1402.0102}
\bibitem{3} T.L. Heath, {\em Diophantus of Alexandria}, Cambridge 1910
\bibitem{4} Walter Wyss, {\em Perfect Parallelograms}, American Math Monthly, 119 (6) (2012), p.513-515
\end{thebibliography}
\end{document}